\documentclass[11pt,letterpaper]{amsart}
\usepackage[english]{babel}
\usepackage[utf8]{inputenc}
\usepackage{lmodern}
\usepackage[T1]{fontenc}
\usepackage{graphicx}
\usepackage{amsmath,amsthm,amssymb}
\usepackage{amsfonts}
\usepackage{enumitem}
\usepackage{stmaryrd}
\usepackage{esint}
\usepackage{soul}
\usepackage{subeqnarray}
\usepackage[all]{xy}
\usepackage{varioref}
\usepackage{tikz-cd}
\usepackage{hyperref}
\usepackage{pgfplots}
\newcommand{\cat}{\text{CAT}}
\newcommand{\ch}{\text{ch}}
\newcommand{\chinf}{\text{ch}(X^\infty)}

\newcommand{\bdinf}{X^\infty}
\newcommand{\Aut}{\text{Aut}}

\DeclareMathOperator{\Out}{\text{Out}}
\DeclareMathOperator{\supp}{\text{supp}}

\DeclareMathOperator{\iso}{Isom}
\newcommand{\R}{\mathbb{R}}
\newcommand{\N}{\mathbb{N}}

\newcommand{\proj}{\text{proj}}
\newcommand{\germ}{\text{germ}}

\newcommand{\res}{\text{Res}}

\newcommand{\bd}{\partial_\infty}
\newcommand{\nui}{\check{\nu}}
\newcommand{\mui}{\check{\mu}}

\newcommand{\Opp}{\text{Opp}}

\newcommand{\CAT}{\mathrm{CAT}}

\theoremstyle{plain}
\newtheorem{thm}{Theorem}[section]

\newtheorem{cor}[thm]{Corollary}
\newtheorem{lem}[thm]{Lemma}
\newtheorem{prop}[thm]{Proposition}
\theoremstyle{definition} 
\newtheorem{Def}[thm]{Definition}
\newtheorem{rem}[thm]{Remark}
\theoremstyle{definition}

\theoremstyle{definition}
\newtheorem*{ackn}{Acknowledgements}
\newtheorem{thmi}{Theorem} 

\newcommand{\newcomment}[4]{%
	\newcounter{#2counter}
	\expandafter\newcommand\csname #1\endcsname[1]{%
		\refstepcounter{#2counter}%
		{\color{#4}(#3\arabic{#2counter})}\marginpar{\scriptsize\raggedright\textbf{\color{#4}(#2 \arabic{#2counter}):} ##1}%
}}

\newcomment{clb}{Corentin}{c}{teal}
\newcomment{jl}{Jean}{jl}{blue}
\newcomment{js}{Jeroen}{j}{red}

\title[A Tits alternative for $\tilde{A}_2$-buildings]{A Tits alternative for $\mathbb{R}$-buildings of type $\tilde{A}_2$}
\author{Corentin Le Bars, Jean L\'ecureux \& Jeroen Schillewaert}
\date{}

\begin{document}
	\maketitle
	\begin{abstract}
		Let $G$ be a group with a non-elementary action on a (not necessarily discrete) $\tilde{A}_2$-buildings. We prove that, given a random walk on $G$, isometries in $G$ are strongly regular hyperbolic with high probability. As a consequence, we prove a Tits alternative for $G$, as well as a local-to-global fixed point result. We also prove that isometries of (not necessarily complete) $\mathbb{R}$-buildings are semi-simple.
	\end{abstract}
	
	\section{Introduction}
	
	The celebrated \emph{Tits alternative} asserts that, in various contexts, the following dichotomy occurs for a finitely generated group $G$: either $G$ is "small" (in a sense depending on the context), or $G$ contains a non-abelian free group. The original statement of Tits \cite{tits72} treats the case of linear groups, where he proves that finitely generated groups not containing free groups are virtually solvable.
	
	Since then, there has been a large amount of work to prove the Tits alternative in various contexts: groups acting on hyperbolic spaces \cite[Theorem 5.3.E]{gromov87}, mapping class groups and $\Out(F_N)$ \cite{bestvina_feighn_handel00} (see also \cite{horbez16}), CAT(0) cube complexes \cite{sageev_wise05}, and many others. A common feature in all these examples is that they share some kind of non-positive curvature, and a long-standing question is to know whether groups acting on CAT(0) spaces satisfy such a Tits alternative. 
	
	An impressive progress in that direction was recently obtained by Osajda and Przytycki, who proved a Tits alternative for finitely generated groups acting on 2-dimensional $\mathrm{CAT}(0)$ complexes with a bound on the order of the cell stabilizers \cite{osajda_przytycki21}. Our first theorem is also another step in that direction: we obtain a Tits alternative for countable groups acting on Euclidean buildings. As a corollary we obtain an result analogous to \cite{osajda_przytycki21} for not necessarily discrete Euclidean buildings of type $\tilde{A}_2$.
	
	We also note that in the case of (simplicial) buildings, although also non-affine ones, the recent preprint \cite{karpinski_osajda_przytycki24} proves a Tits alternative, but under the much stronger assumption that the group acts geometrically on the building, which had been proved by Ballmann and Brin in the case of Euclidean buildings \cite[Theorem F]{ballmann_brin95}. In \cite{stadler24}, S.~Stadler also proves a Tits alternative for geometric actions on Hadamard spaces under the additional assumption that $\Gamma$ satisfies Ballmann's duality condition. 	
	In the following theorem  the word "building" refers to an affine building which is not necessarily simplicial, in the sense of Parreau \cite{parreau00} for example. We say the action of $G$ on $X$ is elementary if either $G$ has a bounded orbit in $X$ or if it has a finite orbit in $\partial_\infty X$. Note that another usual notion of non-elementarity also requires that $G$ does not have an invariant flat. However in our setup an invariant flat will be contained in an apartment, hence its boundary will give rise to a finite orbit in the boundary.

	\begin{thmi}\label{thm Tits alt}
		{A countable group with a non-elementary  action  on a building of type $\tilde{A}_2$ contains a non-abelian free subgroup.}
	\end{thmi}

	\begin{rem}
		By default the action on a building is by automorphism of the building, in the sense of \cite[\S 2.1.13]{rousseau23} for example. The group of automorphisms contains a finite-index subgroup formed by \emph{type-preserving automorphisms}, see \cite[\S2.4.6.1]{rousseau23}. Since our results are clearly implied by the same results for a finite-index subgroup, we will also always assume that our actions are type-preserving. Note that automorphisms of buildings are isometries. Moreover,  by \cite[Lemma 3.5]{schillewaert_struyve_thomas22}, we may assume $X$ is an $\mathbb{R}$-building in which each vertex is special.
	\end{rem}
	
	We also obtain the following corollary to Theorem \ref{thm Tits alt}, which as discussed above is a non-discrete generalization of \cite{osajda_przytycki21}. Following \cite[Section 4]{kapovich24}, we say that the $G$-action is \emph{wandering} if for every point $x \in X$, there exists an neighborhood $U_x$ of $x$ such that the transporter subset 
	$$ (U_x|U_x)_G:= \{g \in G \mid gU_x \cap U_x \neq \emptyset\}$$
	is finite.

	\begin{thmi}\label{elementary}
		Suppose that a finitely generated group $G$ has a wandering action on a building $X$ of type $\tilde{A}_2$ with uniformly bounded point-stabilizers. Then $G$ is virtually cyclic, virtually $\mathbb{Z}^2$ or contains a non-abelian free group.
	\end{thmi}
	
	\begin{rem}
		Kapovich \cite{kapovich24} investigates various notions related to proper discontinuity, of which wandering action is a weakening.
		If the building is simplicial, our only assumption is that point-stabilizers are uniformly bounded.
	\end{rem}

	Affine buildings of dimension at least 3 have been classified by yet another impressive work of Tits: they are all of algebraic origin, and therefore amenable to the original methods of Tits. However, in dimension 2, there are many examples of exotic buildings, not coming from any algebraic construction. 
	This justifies our restriction to dimension 2. Nevertheless, there are two other irreducible types of affine buildings of dimension 2 (namely, $\tilde C_2$ and $\tilde G_2$). Our methods of proof are unable to reach these cases for the moment (see Remark \ref{rem:othertypes} below).
	
	Another problem, very closely related to the Tits alternative (and which is sometimes a first step towards it), is sometimes called the \emph{local to global principle}.
	Let $G$ be a finitely generated group acting by isometries on a metric space $X$. We call the action \emph{locally elliptic} if every element of $G$ fixes a point in $X$.
	A standard question in geometric group theory is the following: when does that imply that the action of $G$ itself is elliptic, \textit{i.e.} is there a global fixed point? The answer is positive in many situations where $X$ is non-positively curved, see for example the discussion in \cite{haettel_osajda22}.
	
	Our second result is a new proof of the following theorem.
	
	\begin{thmi}[{Fixed point}]\label{local-global}
		A finitely generated group wich acts on a building of type $\tilde{A}_2$ such that each element fixes a point has a  global fixed point.
	\end{thmi}

	Once again, similar results were obtained by Norin, Osajda and Przytycki for various subclasses of $\CAT(0)$ triangle complexes \cite{norin_osajda_przytycki22}, which include 2-dimensional simplicial buildings. Furthermore, for not necessarily discrete Euclidean buildings of type $\tilde{A}_2$ and $\tilde{C}_2$ by Schillewaert, Struyve and Thomas \cite{schillewaert_struyve_thomas22} prove the same theorem. The former uses Helly's theorem, Masur's theorem on periodic trajectories in rational billiards, and methods by Ballmann and Brin to find closed geodesics in 2-dimensional locally $\CAT(0)$ complexes, whereas the latter uses intrinsic building theoretic methods using the local spherical buildings in conjunction with general $\CAT(0)$ space techniques including Busemann functions.  Even though Theorem \ref{local-global} is covered by \cite{schillewaert_struyve_thomas22}, the methods of proofs are very different and we have hope that the current approach could be pushed further to a more general context (starting with buildings of type $\tilde G_2$).  Breuillard and Fujiwara prove a quantitative version of this result for Bruhat-Tits buildings based on a Bochi-type inequality \cite[Theorem 7.16]{breuillard_fujiwara21} and ask whether this also holds for the isometry group of arbitrary affine buildings or finite-dimensional $\mathrm{CAT}(0)$ cube complexes. For $\mathrm{CAT}(0)$ square complexes Kar and Sageev \cite{kar_sageev19} gave a positive answer. Another recent result in that direction was obtained very recently by Izeki and Karlsson \cite{izeki_karlsson24}, for general CAT(0) spaces but a more restrictive class of groups, using methods which are more in the spirit of the current paper.
	
	\medskip

	The heart of  the proofs of both theorems is quite different from the classical ones, as it relies on results on random walks on the group $G$. Let $\mu$ be a probability measure on $G$ whose support generates $G$, we consider the associated random walk $(Z_n)$ on $G$: in other words, we pick at random  elements $g_i$ of law $\mu$, independently, and form the product $Z_n = g_1g_2\dots g_n$. 
	
	Let $X$ be a $\mathrm{CAT}(0)$ space, let $g$ be a hyperbolic isometry of X and let $\mathrm{Min}(g) := \{x \in X \mid d(x, g(x)) = |g|\}$,
	where $|g|$ denotes the translation length of $g$. We say that $g$ is a {\em regular} hyperbolic isometry of $X$ if $\mathrm{Min}(g)$ is at bounded Hausdorff distance from a maximal flat of $X$.
	
	A regular hyperbolic isometry $g$ of a building $X$ is called \emph{strongly regular}  if one (and hence all) of its translation axes crosses all the walls of the unique apartment contained in $\mathrm{Min}(g)$. The main interest of strongly regular elements is that they have a weak form of North-South dynamics, which is very useful in producing free groups (see  \cite[Section 2]{caprace_ciobotaru15}).

	Our main result, from which both Theorem \ref{local-global} and \ref{thm Tits alt} follow, is that the proportion of strongly regular hyperbolic isometries in the random walk $\{Z_n\}$ with respect to the action on $X$ goes to 1 almost surely as $n $ goes to infinity. The assumptions on the random walk in the following statement will be defined in Section \ref{sec:statmeas}.

	\begin{thmi}\label{thm prop loxod}
		Let $X$ be a building of type $\tilde{A}_2$, let $G$ be a countable group  and let $G \curvearrowright X $ be a non-elementary action by type-preserving automorphisms. Let $\mu$ be a symmetric and admissible probability measure  with finite second moment, and denote by $(Z_n)$ the associated $\mu$-random walk. Then 
		$$ \mathbb{P}\big( \, Z_n\text{ is a strongly regular hyperbolic isometry}\, \big) \underset{n\to \infty}{\longrightarrow}1. $$
	\end{thmi}

	Theorem \ref{thm prop loxod} is analogous to, for instance, \cite[Theorem 11.5]{fernos_lecureux_matheus18} where it is proven that the proportion of contracting elements for random walks in $\cat(0)$ cube complexes goes to 1, or to \cite[Theorem 1.2]{le-bars24tcl}, where the author shows a similar result for random walks with contracting elements in Hadamard spaces.

	Its proof relies heavily on a recent result obtained by the first author in \cite{le-bars23immeuble}. In particular he proved that, under a second moment assumption on $\mu$, $(Z_n o)$ converges almost surely to a regular (i.e. in the interior of a chamber) point of the boundary, for any base-point $o \in X$ (see Section \ref{subsec:Corentin} for more precise statements).

	\begin{rem}\label{rem:othertypes}
		Theorem \ref{thm prop loxod} hinges on the existence of an opposition involution in type $\tilde{A}_2$ which is non-trivial on every type of vertices. If the opposition involution is trivial on the types of vertices (which is the case for buildings of type $\tilde{C}_2$ or $\tilde{G}_2$), it follows from the results of \cite[Th\'eor\`eme 5.1]{quint02} that there exists groups with a non-elementary action on a building, but with no strongly regular hyperbolic elements. In fact the results of \cite{le-bars23immeuble} are not true in this context.  While we still hope that our methods could be able to find hyperbolic elements and to prove a Tits alternative in the general case, it seems that the necessary adaptations would be highly non-trivial.
	\end{rem}
	
	We conclude the introduction with an extension of a result by Parreau to the metrically incomplete case \cite[Corollaire 4.2]{parreau00}.
	
	\begin{thmi}\label{semi-simple}
		Isometries of $\mathbb{R}$-buildings are semi-simple.
	\end{thmi}

	\begin{ackn}
		We are grateful to Anne Parreau for suggesting a simplification of a previous version of the argument in Lemma \ref{hyperbolic-creation}, and to Koen Struyve for the proof strategy of Theorem \ref{semi-simple}.
		The research of the third author is supported by the New Zealand Marsden Fund through grant UOA-2122.
	\end{ackn}
	
	\section{Affine buildings}
	\label{section intro immeuble}
	
	In this section, we define non-discrete affine buildings, and briefly discuss some of their basic properties, we refer the reader to \cite{parreau00} or \cite{rousseau23}
	for more information.

	\subsection{Non-discrete affine buildings}\label{section intro affine buildings}
	
	Let $(W,V)$ be an affine reflection system, that is $W = W_0 \ltimes T$, where $W_0$ is a finite reflection group and $T$ is a translation group on $V=\mathbb{E}^n$. We denote by Let $\mathfrak{a}^{+}$ some Weyl chamber for $W$, and by $\mathfrak{a}^{++}$ its interior. Let $X$ be a set, and let $\mathcal{A}$ be a collection of injective charts of $V$ into $X$, which we call an \emph{atlas}. Each such injection is called a \emph{chart}, or \emph{marked apartment}, and the image $A$ of $V$ by an injection is called an \emph{apartment}. We say that $(X, \mathcal{A})$ is an affine building modeled after $(W, V)$ if the following axioms are verified. 
	\begin{enumerate}[label=(A\arabic*)]
		\item The atlas $\mathcal{A}$ is invariant by pre-composition with $W$. \label{A1}
		\item Given two charts $f, f' : V \to X$ with $f(V) \cap f'(V) \neq \emptyset$, then $U:=f^{-1} (f'(V))$ is a closed convex subset of $V$, and there exists $w \in W$ such that $f|_U = f' \circ w |_{U}$. \label{A2}
		\item For any pair of points $x, y \in X$, there is an apartment containing both. \label{A3} 
	\end{enumerate}
	Axioms \ref{A1}-\ref{A3} imply the existence of a well-defined distance function $d : X \times X \to \R_+ $, such that the distance between any two points is the $d_V$-distance between their pre-image under any chart containing both. The metric space $(X,d)$ is then a CAT(0) space. Every automorphism of $X$ induces an isometry of $(X,d)$.
	A \emph{Weyl chamber} (or \emph{sector}) in $X$ is the image of an affine Weyl chamber under some chart $f \in \mathcal{A}$. 
	\begin{enumerate}[label=(A\arabic*), resume]
		\item Given two Weyl chambers $S_1, S_2$ in $X$, there exist sub-Weyl chambers $S_1'\subseteq S_1, S_2' \subseteq S_2$ such that $S_1'$ and $S_1'$ are contained in the same apartment. \label{A4} 
		\item For any apartment $A$ and $x \in X$, there exists a retraction $\rho_{A, x}: X \to A$ such that $\rho_{A, x} $ does not increase distance and $\rho^{-1}_{A, x} (x ) = \{x\}$. \label{A5} 
	\end{enumerate}
	
	If the affine reflection group is not discrete, we say that the building $(X,\mathcal{A})$ modelled after $(W, V)$ is \emph{non-discrete}. We will assume that the system of apartments $\mathcal{A}$ is maximal. We say that $X$ is of type $\tilde{A}_2$ if $W_0$ is a spherical Coxeter group of type $A_2$. 
	
	By the axioms \ref{A1}-\ref{A3}, there exists a marked apartment $f$ sending the fundamental closed Weyl chamber $\mathfrak{a}^+$ to a Weyl chamber in $X$ based at $x $ and containing $y$. The \emph{type} $\theta(x,y)$ of the Euclidean segment $[x,y]$ is the unique vector in $\mathfrak{a}^+$ such that $y = f(\theta(x,y))$. The group $G$ is type-preserving by assumption: this implies that for any $g\in G$, $\theta(gx,gy)=\theta(x,y)$.

	\subsection{Spherical building at infinity}
	
	As a $\cat$(0) space, $X$ has a \emph{visual} bordification, given by equivalence classes of rays, two rays being equivalent if they are at finite Hausdorff $d$-distance, denoted by $\overline{X} = X \cup \bd X$. The visual boundary $\bd X$ can be endowed with a natural topology and isometries of $X$ extend to homeomorphisms on the boundary.  Moreover, when $X$ is separable, the topology on $\bd X$ is metrizable. 
	Analogous to defining ends of trees one can equip $\bd X$ with the structure of a spherical building \cite[Propriété 1.7]{parreau00}.
	
	\subsection{The cone topology on $\overline{X}$.} 
	
	The visual boundary is a classical notion for complete CAT(0) spaces, see e.g \cite{bridson_haefliger99}. For non-complete spaces it is usually not well-defined (for example, it might depend on a basepoint). However for Euclidean buildings, even non-complete, it turns out that the construction still works well, see \cite{rousseau23}. We recall below the definition.
	
	Fix a point $x_0\in X$ and for $x\notin B(x_0,r)$ define $p_r(x)$ to be the unique point on $[x_0,x]$ at distance $r$ from $x_0$. Consider the inverse system of projection maps $p_r\mid_{\bar{B}(x_0,r)}$.
	We consider $\varprojlim \bar{B}(x_0,r)$ with the inverse limit topology. A point in this space is given by a map $c:[0,\infty)\to X$ such that if $r'\geq r$, then $p_r(c(r')) = c(r)$ and the inverse limit topology
	coincides with the topology of uniform convergence on compact subsets. The \emph{cone topology} is the topology $\mathcal{T}(x_0)$ such that the natural bijection $\phi(x_0): \bar{X}\to \varprojlim \bar{B}(x_0,r)$ is a homeomorphism. This topology is independent of the basepoint $x_0$. When we equip $\partial_{\infty} X$ with this topology, we call it the {\emph{visual boundary}.
		
		Note that for $x$ a vertex, and $F^\infty$ a facet at infinity, there is an affine facet based at $x$ in the equivalence class of $F^\infty$, which we denote by $Q(x, F^\infty)$. 
		
		A detailed account on how to give a topology on the set $X \cup \chinf$, and even on $X \cup X^\infty_\tau$, where $X^\infty_\tau$ represents the set of simplices of type $\tau$ of the spherical building at infinity, is given in \cite{rousseau23}. We will not give many details on the topology, but we summarize here a few results. 
		
		A basis of open neighborhoods in $\ch (\bdinf)$ is given by 
		\begin{eqnarray}\label{eq basis chinf}
			U_x(y) := \{ C \in \ch(\bdinf) \, | \, y \in Q(x, C)\} \subseteq\chinf, 
		\end{eqnarray}
		for $x,y \in X$. With this topology, the set of chambers at infinity $\ch (\bdinf)$ is a totally disconnected space.
		
		The following proposition summarizes some properties of this bordification.  The set of \emph{regular points} $\bd^{\mathrm{reg}} X$ is defined as the set of boundary points $\xi \in\bd X$ that are strictly supported on a chamber, denoted $C_\xi$.
		
		\begin{prop}[{\cite[\S 3.2]{rousseau23}}]
			Let $X$ be any Euclidean building. Then there is a topology on $X \cup \chinf$ for which a basis of open sets of the chambers at infinity is given by the sets \eqref{eq basis chinf}. It agrees with the $\cat$(0) topology on $X$. This topology is first-countable and Hausdorff.
			The map $\bd^{\mathrm{reg}} X	\to \chinf$ defined by $\xi\mapsto C_\xi$ is a homeomorphism.
		\end{prop}
		
		Two chambers in a spherical building are said to be \emph{opposite} if the gallery distance between them is maximal, see for instance \cite[Chapter 5]{weiss03}. In this case, there is a unique apartment joining them. 
		
		\subsection{Retractions}
		
		Let $X$ be an affine building, and let $C\in \chinf$. If $A$ is an apartment such that $C\in \partial A$ (equivalently, there exists a Weyl chamber $S$ contained in $A$ representing $C$). Then there exists a unique retraction map $\rho_{A, C}: X \to A$ such that $\rho_{A, C}$ preserves the distance on any apartment containing a chamber representing $C$, see \cite[Proposition 1.20]{parreau00}. Moreover, $\rho_{A, C}$ does not increase distances. We call $\rho_{A,C}$ the \emph{canonical retraction} of $X$ on $A$ based at the chamber at infinity $C \in \chinf$. In particular, for any apartment $A'$ containing $C$ in its boundary, ${\rho_{A,C}}_{\mid A'} : A' \to A$ is an isomorphism fixing $A \cap A'$ pointwise.

		\subsection{Regularity} 
		
		Again, let us denote by $\mathfrak{a}^{+}$ the fundamental Weyl chamber, and denote by $\mathfrak{a}^{++}$ its interior. Denote by $w_0$ the long element of the finite Weyl group $W_0$ associated with the affine reflection group $W$. For $\lambda \in \mathfrak{a}$, the opposition involution $j : \mathfrak{a}^{+} \to \mathfrak{a}^{+}$ is defined by 
		$$ j (\lambda) = w_0(-\lambda).$$ 
		Let $x, y  \in X$ be points of the building $X$. By Axiom (A3), there is a marked apartment $f \in \mathcal{A}$ containing both, and by Axiom (A2) we can assume that $f(0)= x $ and $f^{-1}(y) \in \mathfrak{a}^{+}$. We say that the segment $[x,y]$ is \emph{regular} if the type $\theta(x,y)$ is regular, i.e. $\theta(x,y) = f^{-1}(y) \in \mathfrak{a}^{++}$. 
		If $[x,y]$ is regular, then both $\theta(x,y)$ and $j(\theta(y,x))$ are regular. Notice that in this case, $f(\mathfrak{a}^+)$ is a Weyl chamber of the building $X$, contained in $A$ and containing $y$. Denote by $C_y$ the chamber at infinity that it represents. Similarly, $f( \theta(x,y)- \mathfrak{a}^+)$ is a Weyl chamber containing $x$ representing a chamber $C_x \in \chinf$. These chambers are opposite in the spherical building at infinity. In this case, there exists a unique apartment $A$ joining them, i.e. such that both chambers belong to $A^\infty$. 
		
		The following proposition is from \cite[Proposition 2.10]{caprace_ciobotaru15}, see also \cite[Proposition 6.1]{ciobotaru_muhlerr_rousseau_20} for the non-discrete case. It roughly states that strongly regular hyperbolic elements satisfy a weak version of North-South dynamics, similar to loxodromic isometries in hyperbolic spaces or contracting isometries in $\cat$(0) spaces.
		
		\begin{prop}
			Let $g \in \iso(X)$ be a type preserving strongly regular hyperbolic element, with unique translation apartment $A$. Let $C^-, C^+ \in \chinf$ be the repelling (resp. attracting) chamber at infinity for $g$. Then for every $C \in \chinf$, the limit $\lim g^n (C) $ exists and coincides with $\rho_{A, C^-} (C)$, where $\rho_{A, C^-}$ is the retraction onto $A$ centered at $C^-$. 
		\end{prop}

		\subsection{Residue building based at a vertex}
		Non-discrete buildings still possess a good notion of (local) alcoves and faces. 
		Let $F$ and $F'$ be two facets based at a vertex $o \in X$. We say that $F$ and $F'$ have the same \emph{germ at $o$} if their intersection is an open neighborhood of $o$ in both $F$ and $F'$. The set of all germs at $o$ that are not reduced to $\{o\}$ can be given the structure of a simplicial spherical building, called the \emph{residue building} based at $o$ \cite[Corollary 1.11]{parreau00}. Simplices of maximal dimension of $\Sigma_o X$ are called (local) \emph{alcoves}. 
		
		For every vertex $x \in X$, we have a canonical morphism of simplicial complexes 
		\begin{eqnarray}
			\Sigma_x : \bdinf \to \Sigma_x X \label{morph spherical residue}
		\end{eqnarray}
		sending any facet at infinity $F^\infty $ to $\germ_x(Q(x,F^\infty))$.

		Fix $o \in X$. For every $ x \in X$ such that $[o,x]$ is regular, there exists $C \in \chinf $ such that $x $ belongs to the interior of $ Q(o, C)$, we can associate a unique local alcove in $\Sigma_o X$ defined by $\Sigma_o(x)=\germ_o (Q(o, C))$. The following definition appears in \cite{caprace_lecureux11}. 
		
		\begin{Def}
			We say that a sequence $(x_n)$ converges to the chamber at infinity $C \in \chinf$ \emph{in the combinatorial sense} if for every $o \in X$, there exists $n_0$ such that for all $n \geq n_0$, the projection $\Sigma_o(x_n)$ of $x_n$ on the residue building $\Sigma_o X$ is the chamber $\Sigma_o (C)=\germ_o(Q(o,C))$. 
		\end{Def}
		
		For more on this notion, we refer to \cite{caprace_lecureux11} (in the discrete case) or to \cite[\S 3.2]{rousseau23} for non-discrete affine buildings. 
		
		\subsection{Panel trees}\label{section panel trees}
		Recall that a boundary point $v \in  \bd X$ is an equivalence class of rays, two rays $r_1$ and $r_2$ being equivalent, for which we write $r_1 \sim r_2$, if they contain subrays that lie in a common apartment and are parallel in this apartment.  We will say that two geodesic rays $r_1$ and $r_2$ are \textit{strongly asymptotic}, and write $r_1 \simeq r_2$, if their intersection contains a geodesic ray.  For two equivalent geodesic rays $r_1$ and $r_2$ that represent the boundary point $v \in \bd X$, we define their distance to be:
		
		\begin{equation*}
			d_v (r_1, r_2) := \underset{s}{\inf} \underset{t \rightarrow \infty}{\lim} d(r_1(t+s), r_2(t)). 
		\end{equation*}
		Note that it defines a pseudo-distance \cite[Section 5.3]{caprace_lecureux11}, and that two strongly asymptotic rays  $r_1 $ and $r_2$ satisfy $d_v(r_1, r_2) =0$. This pseudo-distance does not depend on the $\simeq $-strongly asymptotic classes of rays among rays from the same $\sim$-equivalence class: on these $\simeq$-classes, it becomes a distance.
		
		For an affine building of dimension 2 and $v$ a vertex at infinity, the metric space $(T_v, d_v)$ of asymptotic classes of rays in the class of a vertex at infinity is an $\R$-tree (an affine $\tilde{A}_1$-building) called the \emph{panel tree} at $v$. The branch points of this tree correspond to thick walls of $X$ \cite[\S 4.26]{kramer_weiss14}.
		Note that given a vertex $v \in \bdinf$, there is a well-defined and continuous application defined by 
		\begin{eqnarray}
			\pi_v : X &\longrightarrow& T_v \nonumber \\
			x &\longmapsto & [Q(x, v)]\nonumber,
		\end{eqnarray}
		where $[Q(x, v)]$ is the class (for the strongly asymptotic relation) of the geodesic ray based at $x$, in the direction of $v$. 
		
		By \cite[Proposition 4]{tits86} there is a canonical $\mathrm{Aut}(X)_v$-equivariant bijection
		between the set of ends of the panel tree $T_v$ and the set $\mathrm{Ch}(v)$ of chambers of $\partial X_\infty$ containing $v$.
		
		\subsection{Reduction: separability}
		
		The following proposition will be used to reduce to the case when $X$ is a separable building.
		
		\begin{prop}\label{prop:separable}
			If $\Gamma$ is a countable group acting on an affine building $X$, then there is a subbuilding $Y\subset X$, of the same type, which is separable and $\Gamma$-invariant.
		\end{prop}
		
		In the sequel we fix such an $X$ and $\Gamma$.
		
		\begin{lem}\label{sub-inv}
			There exists a $\Gamma$-invariant countable sub-building $Z$ of the same type as $X^\infty$.
		\end{lem}
		
		\begin{proof}
			Let $S_0=A$ be an apartment in $\partial_\infty X$, and let $\Gamma(S_0)$ be its $\Gamma$-orbit. If $\Gamma(S_0)$ is a sub-building we are done. Suppose otherwise. Given two simplices of $\Gamma(S_0)$ for which there is no apartment containing them yet add such an apartment. Call the newly obtained set $S_1$. Now repeat the above construction process with $S_1$ and so on. Then $\bigcup_{i=0}^\infty \Gamma(S_i)$ is countable and $\Gamma$-invariant by construction. Moreover it is a subbuilding of the same type as $X^\infty$.		 
		\end{proof}
		
		\begin{proof}[Proof of Proposition \ref{prop:separable}]
			Let $Z$ be a countable invariant subbuilding of $\partial_\infty X$ which exists by Lemma \ref{sub-inv}. Let $Y$ be the union of all apartments whose boundary is entirely contained in $Z$. By construction $Y$ is $\Gamma$-invariant. Since between two opposite chambers of $Z$ there is a unique apartment, $Y$ is a countable union of apartments, hence is separable.
			
			We have to prove that $Y$ is indeed a sub-building of $X$. Let $x\in Y$ be a vertex, and $C$ be a chamber of $Z$. We first claim that there exists an apartment of $Y$ containing the sector $Q(x,C)$. Indeed, since $x\in Y$ there is an apartment $A$ containing $x$ whose boundary is in $Z$. In particular $\proj_x(A)$ is an apartment in $\Sigma_x X$ \cite[Proposition 1.14]{parreau00}. Hence $\Sigma_x X$ contains a local alcove which is opposite $\proj_x(C)$. Let $C'$ be the chamber in the boundary of $A$ which projects to this local alcove. Then $C$ and $C'$ are opposite at $x$, hence there is an apartment containing $x$, $C$ and $C'$. Since $C$ and $C'$ are opposite the apartment of $\partial_\infty X$ containing them is entirely contained in $Z$, and therefore the apartment whose boundary contains both is in $Y$. Since this apartment contains $Q(x,C)$, we proved the claim. 
			
			Now let $x,y\in Y$, we want to prove that $x$ and $y$ are contained in an apartment of $Y$. 
			Then there exists a chamber $C$ in $Z$ such that $\proj_y(C)$ is opposite $\proj_y(x)$ in $\Sigma_x X$. Thus there exists an apartment $A'$ of $X$ containing $Q(y,C)$ and the convex hull of $x$ and $y$. In $A'$ we see that $y\in Q(x,C)$. By the claim above there exists an apartment of $Y$ containing $Q(x,C)$, and this apartment contains $x$ and $y$.
		\end{proof}
		
		\subsection{Reduction: completeness}
		
		We let $X$ be a possibly non-discrete or even non-complete Euclidean building. 
		Let $\bar X$ be the completion of $X$.
		
		\begin{prop}\label{prop vis bd completion}
			The visual boundaries $\bd X $ and $\bd \bar X $ are equal. 
		\end{prop}
		
		\begin{proof}
			Note first that by  \cite{kleiner_leeb97} (see also \cite[Lemma 4.4]{struyve11}) we can, using ultrapowers, embed $X$ into a complete Euclidean building, which we denote $\hat X$.
			
			We first claim that $\bd \bar X$ is the set of all limit points in $\partial_\infty \hat X$ of sequences of points of $X$. Indeed, one inclusion is clear. For the converse, recall that the set $\bd \bar X$ can be defined as the set of geodesic rays in $\bar X$ starting from the point $o \in X $ (the topology on $\bd \bar X $ does not depend upon the choice of the basepoint $o$, see for instance \cite[\S 3.2]{rousseau23}). If $\xi \in \bd \bar X$, then $\xi$ is a limit of points on a geodesic ray of $\bar X$. Each of these points is arbitrarily close to a point of $X$. So $\xi$ is a limit point (for the topology on $\hat X \cup \bd \hat X$) of some sequence of points of $X$. 
			
			Now assume that $(x_n)_{n\in\N}$ is a sequence of points of $X$ converging to $\xi\in \bd \bar X$, with $\xi$ in the interior of a chamber at infinity $C \in \ch(\hat{X}^\infty)$ (of the spherical building at infinity $\hat{X}^\infty$ of the Euclidean building $\hat X$). Let $o\in X$ be a vertex. Then by definition the geodesic segments $[o,x_n]$ converge (uniformly on every compact) to the geodesic ray $[o,\xi)$ which is in the interior of the sector $Q(o,C)$. In particular, for $n$ large enough, the germ of $[o,x_n]$ is in the interior of the local alcove $\proj_o(C)$. It follows that for every $o \in X$ we have $\Sigma_o(C)$ is a local alcove of the residue building $\Sigma_oX$. We claim that this implies that $C \in \ch(X^\infty)$. 
			
			Indeed, let $(x_n)$ be as before a sequence of points in $X$ that converge to $\xi$. For any $o \in X$, the geodesic ray $[o, \xi)$ is contained in the interior of the Weyl chamber $Q(o, C)\subseteq \hat X$. In particular, the sequence $(x_n)$ is asymptotically uniformly $\sigma_{\text{mod}}$-regular in the sense of \cite[Definition 3.1]{kapovich_leeb_porti18morse}. Fix $o \in X$. Consider the geodesic ray $c :[o, \infty) \to \bar X$ starting from $o$ representing $\xi$. For any $n$, pick $x_n \in X$ such that $d(x_n, c(n)) \leq \varepsilon$, for $\epsilon $ fixed (where we still denote by $d$ the metric on $\hat X$). Let $C_n \in \ch(X^\infty)$ be chambers at infinity of $X$ such that for all $n$,  
			$$ C_n \in \{D \in \ch(X^\infty) \mid x_n \in Q(o, D)\}.$$
			Such chambers exist because of axiom \ref{A3} for the Euclidean building $X$, and because $[o, x_n]$ is a regular segment in $X$. Now by convexity of the distance in $\cat$(0) spaces, we have that for every $n \leq m \in \N$, $x_n $ is contained in the $2\varepsilon$ neighborhood $N_{2\epsilon}(Q(o, C_m))$ of $Q(o, C_m)$. Therefore, we can apply \cite[Lemma 3.79]{kapovich_leeb_porti18morse}: there exists $C' \in \ch(X^\infty)$ such that $(x_n)$ converges to $C'$ conically. At the same time, since $(x_n)$ is an asymptotically uniformly regular sequence contained in a tubular neighborhood of $Q(o,C)$, \cite[Lemma 3.76]{kapovich_leeb_porti18morse} ensures that $(x_n)$ converges to $C$ conically. By uniqueness of the limit, we get that $C = C' \in \ch(X^\infty)$. 
			
			Putting everything together we get that $\bd \bar X$ does not contain any regular point (of the spherical building $\hat{X}^\infty$) which is not already in $\partial_\infty X$. Now we know that $\bd \bar X$ is the boundary of a complete CAT(0) space. In particular it is $\pi$-geodesic for the angular metric \cite[Theorem II.9.13]{bridson_haefliger99}. If $\xi'\in \bd\bar X$ is not regular, then there exists a chamber $C\subset\partial_\infty X$ which contains a regular point $\xi$ with $\angle(\xi,\xi')<\pi$ (where $\angle$ is the angular metric on $\partial_\infty \hat X$). Let $C_0$ be the projection of $C$ on the facet containing $\xi$ (in $\partial_\infty \hat X$). By $\pi$-geodesicity of $\bd\bar X$,  the $\angle$-geodesic ray from $\xi'$ to some point in $C$ is in $\bd \bar X$.  Thus this geodesic 
			ray intersects the chamber $C_0$ of $\partial_\infty \hat X$. In particular $C_0$ contains a point $\xi''\in \bar \hat X$, and by the previous argument it follows that $C_0\in \bd X$. But $\xi'$ is contained in the boundary of $C_0$, and since $\partial_\infty X$ is a spherical building we get that $\xi'\in \partial_\infty X$.
		\end{proof}
		
		\begin{cor}\label{cor:nonelementaire}
			If $G $ is a group with a non-elementary isometric action  on $X$, then the $G$-action on $\bar X$ is also non-elementary.
		\end{cor}
		
		\begin{proof}
			If the action is elementary on $X$, then there is a bounded orbit in $X$ (hence in $\bar X$) or a finite orbit in $\bd X=\bd\bar X$, so that the action on $\bar X$ is non-elementary.
			
			If the $G$-action on $\bar X$ is elementary, then either there is a finite orbit in $\bd X=\bd\bar X$ and the action on $X$ is elementary, or 
			there is a bounded orbit (hence a fixed point) in $\bar X$, and since $G$ acts by isometries there is a bounded orbit in $X$, so that the action is again elementary.
		\end{proof}

		We can also deduce that isometries are semi-simple from Proposition \ref{prop vis bd completion}. 
		
		\begin{proof}[Proof of Theorem \ref{semi-simple}]
			Let $g$ be an isometry of the $\R$-building $X$. As in the proof of Proposition \ref{prop vis bd completion} we start by embedding $X$ in a metrically complete $X'$. Since $X'$ is complete the action of $g$ on $X'$ is either elliptic or hyperbolic \cite[Corollaire 4.2]{parreau00}. 
			
			Assume first that the action of $g$ on $X'$ is elliptic, that is, fixes a point of $X'$. Then all $g$-orbits on~$X$ are bounded.  Hence, as $G=\langle g \rangle$ is finitely generated, it fixes a point of $X$ by \cite[Main Theorem 1]{struyve11} and is thus elliptic.
			
			Hence we may assume the action of $g$ on $X'$ is hyperbolic. By \cite[II.6.8(1)]{bridson_haefliger99} there exists a translation axis $L$ whose endpoints $x^+$ (the attracting fixed point) and $x^{-}$ (the repelling fixed point) on the boundary $\partial X'$ of $X'$ are fixed by $g$. 
			
			As $\bar{X}$ is a closed convex complete subset of $X'$, there exists a projection $\pi : X' \to X$ such that $\forall x \in X'$, $\pi(x)$ is uniquely defined and $d(x, \pi(x))= d(x, \bar{X}) = \inf_{y \in \bar{X}} d(x, y) $ \cite[Lemma 2.4]{bridson_haefliger99}. Let $x\in X'$.
			Since $g$ is an isometry of $X$, it leaves the set $\bar{X}$ invariant. In particular, we have that $d(x,\pi(x))=d(x,\bar{X}) = d(g(x),\bar{X})$. Since $d(x,\pi(x))=d(g(x),g(\pi(x))$ we have by uniqueness of the projection that $\pi(g(x)) = g(\pi(x))$. 
			
			Now if $x \in \mathrm{Min}(g)$, we have that $d(\pi(x), g(\pi(x))) = d(\pi(x), \pi(gx)) \leq d(x, gx) $ since $\pi$ is distance non-increasing. Therefore, $\pi(x) \in \mathrm{Min}(g)$. In particular, the $\langle g \rangle$-translates of the geodesic segment $[x, gx]$ define an axis. This axis is contained in $\bar{X}$ because $\bar{X}$ is convex and $g$-invariant. 
			Hence we may assume that the axis $L$ belongs to $\bar{X}$. 
			
			Since $\partial_\infty \bar{X}=\partial_\infty X$ by Proposition \ref{prop vis bd completion} we obtained that $\{x^-,\:x^+\} \subset \partial X$.
			
			Consider the union $Y$ (resp. $Y'$) of all geodesic lines  in $X$ (respectively $X'$) with endpoints equal to $\{x^-,\:x^+\}$. By \cite[Proposition 4.8.1]{kleiner_leeb97}, the spaces $Y$ and $Y'$ are sub-building of $X$ and $X'$ respectively (in fact they are products of a flat and a thick building). Moreover, $\mathrm{rank}(Y) = \mathrm{rank}(X)-1$ and $\mathrm{rank}(Y') = \mathrm{rank}(X')-1$. Note that when we remove the flat factor we get the more standard construction of a transverse building, described for example in 
			\cite[Theorem 3.3.14.1]{rousseau23}. Also note that $Y$ is by construction a sub-building of $Y'$. 
			
			The induced action of $g$ on $Y'$ has a fixed point (corresponding to $L$) and hence the action on $Y'$ and thus also $Y$ has bounded orbits. By \cite[Main Result 1]{struyve11}  $g$ has a fixed point $p$ on $Y$. Let $M$ be the axis in $X$ corresponding to $p$. By the Flat Strip Theorem \cite[II.2.13]{bridson_haefliger99} the axes $L$ and $M$ bound a flat strip. Since $g$ acts as a translation on both $L$ and $M$ we obtain $M\subset \mathrm{Min}(g)$, thus the action of $g$ on $X$ is hyperbolic. 
		\end{proof}

		\section{Stationary measures}\label{sec:statmeas}
		
		The notion of stationary measure is a fundamental tool in the study of random walks on groups. After recalling the standard general background we explain here the main results from \cite{le-bars23immeuble}.
		
		\subsection{Generalities}

		Let us first recall the general background on random walks and notations that we will use. Let $G$ be a discrete countable group, and let $\mu$ be a probability measure on $G$. Let $(\Omega, \mathbb{P}) $ be the probability space $(G^{\mathbb{N}}, \mu^{\otimes \mathbb{N}})$, with the product $\sigma$-algebra. The space $\Omega$ is called the space of (forward) increments. The application 
		\begin{equation*}
			(\omega,n) \in  \Omega \times\mathbb{N} \mapsto Z_n(\omega) = \omega_1 \omega_2 \dots \omega_n,
		\end{equation*}
		defines the random walk on $G$ generated by the measure $\mu$.

		Let $(X,d)$ be a metric space endowed with its Borel $\sigma$-algebra and let $\phi : G \to \iso(X)$ be an isometric $G$-action. The data $(\Omega, \mathbb{P}, S, f, \phi)$ defines a measurable random dynamical system 
		\begin{eqnarray}
			& & \Omega \times \mathbb{N} \times X \to X \nonumber \\
			& & (\omega, n, o) \mapsto \phi(Z_n(\omega))  o. \nonumber
		\end{eqnarray}
		If the context is clear, we omit the writing of ``$\phi$'' and we just denote the random walk on $X$ based at $o$ by $(Z_n o)$. 
		If $G$ acts on a Borel space $(Z,\nu)$, the \emph{convolution} $\mu\ast \nu$ is the pushforward of $\mu \otimes \nu$ by the action map $G\times Z\to Z$. In other words, 	if $f $ is a bounded measurable function on $Z$, then
		\begin{equation*}
			\int_Z f(x)d(\mu \ast \nu)(x) = \int_G\int_Z f(g \cdot x ) d\mu(g) d\nu (x).
		\end{equation*}

		We will denote by $\mu_m = \mu^{\ast m}$ the $m$-th convolution power of $\mu$, where $G$ acts on itself by left translation $(g,h) \mapsto  gh$. We want our random walk to visit the whole group, which is why we assume that $\mu$ is  \emph{admissible} measures, that is, its support $\supp(\mu)$ generates $G$ as a semigroup.

		Furstenberg introduced the notion of stationary measures, i.e. measures that are ``invariant'' under the dynamical system. 
		If $(Y, \mathcal{Y})$ is a standard Borel $G$-space, a probability measure $\nu$ on $Y$ is said to be \emph{$\mu$-stationary} if it satisfies 
		$$ \nu = \mu \ast \nu.$$
		
		\subsection{Stationary measures on $\tilde{A}_2$-buildings}\label{subsec:Corentin}
		
		The main results of \cite{le-bars23immeuble} can be summarized as follows. Here, we let $X$ be a separable complete building of type $\tilde{A}_2$, $G$ a  countable group and $G \curvearrowright X$ a non-elementary action on $X$ by isometries. Let $\mu$ be a symmetric and admissible probability measure with finite second moment on $G$.
		
		\begin{thm}[{\cite[Theorem 1.4]{le-bars23immeuble}}]\label{thm:LB24}
			There exists a unique $\mu$-stationary measure $\nu$ on $X\cup \ch(\bdinf)$, and the measure $\nu$ gives full mass to $\ch(\bdinf)$. Furthermore, for $\mathbb P$-almost every $\omega\in\Omega$ and any basepoint $o \in X$, the trajectory $(Z_n(\omega) o)_{n\in\N}$ converges to a regular point $\xi(\omega)$ of the visual boundary, which belongs to the interior of some chamber at infinity $C_\omega$. The law of $C_\omega$ is exactly the measure $\nu$.
		\end{thm}

		We will also use the following fact about the stationary measure.
		
		\begin{thm}[{\cite[Proposition 1.3]{le-bars23immeuble}}]\label{thm:asopposite}
			Let $\nu$ be the unique $\mu$-stationary measure on $\ch(\bdinf)$ given by Theorem \ref{thm:LB24}. Then $\nu \otimes \nu$-almost every pair of chambers in $\ch(\bdinf)$ are opposite. 
		\end{thm}
		Finally, the following fact shows that the convergence of the random walk takes place in the combinatorial sense too. 
		
		\begin{prop}[{\cite[Corollary 7.6]{le-bars23immeuble}}]\label{prop cv comb}
			The random walk $(Z_n (\omega)o )$ converges almost surely to the chamber $C_\omega \in \chinf$ in the combinatorial sense, where $C_\omega$ is given by Theorem \ref{thm:LB24}.
		\end{prop}
		
		\begin{rem}\label{rem:subbuilding}
			In fact, the methods of \cite{le-bars23immeuble} extend to the slightly more general case where $X$ is a separable Hadamard $G$-invariant subspace of a building $Z$ of type $\tilde A_2$, such that the boundary of $X$ is a spherical sub-building of $Z$, see \cite[Theorem 1.7]{le-bars23immeuble}.
		\end{rem}

		\section{Hyperbolic elements}\label{section hyp eleme A2}
		
		In this section we assume that either $X$ is a metrically complete separable building of type $\tilde A_2$, or by Remark \ref{rem:subbuilding}
		that $X$ is a convex, separable, complete subspace of a building of type $\tilde A_2$ whose boundary is a spherical building.
		
		\subsection{Proportion of strongly regular hyperbolic elements}
		
		In the course of the proof, we will use the following lemma, which gives an easy criterion for an isometry to be (strongly regular) hyperbolic in an affine building. In the case when $X$ is a convex, separable, complete subspace of a building  $\hat X$ of type $\tilde A_2$ whose boundary is a spherical building, by a \emph{strongly regular hyperbolic isometry} we mean a strongly regular hyperbolic isometry of $\hat X$, all of whose axes lie in $X$.
		
		\begin{lem}\label{hyperbolic-creation}
			Let $g $ be a type-preserving automorphism of $X$, and assume that there exists a vertex $o \in X_s$ such that $go\neq o$ and  the segments $[g^{-1}o,o]$ and $[o,g(o)]$ are contained in the opposite sectors $S_1,S_2$ respectively, both of which are based at $o$. Then $g$ is a hyperbolic isometry. If moreover $\theta(o, go)$ is regular, then $g$ is a strongly regular hyperbolic isometry. Its translation apartment is the unique apartment containing $S_1$ and $S_2$. 
		\end{lem}
		
		\begin{proof}
			Since $S_1 $ and $S_2$ are opposite sectors based at $o$, there exists a unique apartment $A$ containing them \cite[Proposition 1.12]{parreau00}. By Parreau's axioms (A1) and (A2), we can identify $A$ with $V$ by an isometry (the inverse of a marked apartment) $f_{A,S_2} : A \to V$ sending $S_2 $ to the fundamental Weyl chamber $\mathfrak{a}^{++}$. Since $g$ is type-preserving, the types of $\theta(g^{-1}o, o)$ $\theta(o,go)$ are the same. Since $\theta(g^{-1}o , o)= j(\theta(o,g^{-1}o))$, where $j$ is the involution on the fundamental Weyl chamber, the distinct points $f_{A,S_2}(g^{-1}o), f_{A,S_2}(o)$ and $f_{A,S_2}(go)$ must be aligned. If we repeatedly apply $g$ and $g^{-1}$ to this segment, we see that all the iterates of $o$ by $g$ are aligned, and therefore that $g$ is hyperbolic. The direction of its translation axes is given by $\theta(o, go)$. Therefore, if $\theta(o,go)$ is in the interior of the fundamental Weyl chamber, the translation axes are regular and by \cite[Lemma 2.3]{caprace_ciobotaru15}, $g$ is strongly regular hyperbolic. Note that  \cite[Lemma 2.3]{caprace_ciobotaru15} is only stated for simplicial buildings, but the proof extends verbatim to our context.
			
			In case when $X$ is a subset of a building $\hat X$,  a priori the axes of $g$ could be in $X$, but its attractive and repelling endpoints are in $\bd X$, and by convexity the unique apartment formed by the union of its axes must be contained in $X$, so $Z_n$ is strongly regular with axis in $X$.
		\end{proof}
		
		The following corollary is immediate from the definition. 
		\begin{lem}\label{lem control attracting chamber}
			Let $g \in \iso(X) $ be a strongly regular hyperbolic isometry. Let $\gamma$ be an axis of $g$ passing through a vertex $o \in X$ such that $\gamma (+\infty) $ belongs to the interior of the attracting chamber at infinity $C^+ \in \chinf$, and assume that $\gamma(0)=o$. Then for all $t >0$, $C^+ \in U_o(\gamma(t)) $. 
		\end{lem}
		\begin{proof}
			Let $A$ be the unique translation apartment of $g$. Then $\gamma$ is contained in $A$, and so is $o$. Now for all $t >0$, $\gamma(t) \in Q(o, C^+)$ and the lemma follows. 
		\end{proof}
		\begin{cor}\label{cor control attr ch}
			Let $g $ be a type-preserving isometry of $X$, and assume that there exists a vertex $o \in X_s$ such that $go\neq o$ and  the segments $[g^{-1}o,o]$ and $[o,g(o)]$ are contained in the opposite sectors $S_1,S_2$ respectively, both of which are based at $o$. Assume that $go $ belongs to the interior of $S_1$. Then $g$ is a strongly regular hyperbolic isometry and if we denote by $C^+$ its attracting chamber at infinity, then $C^+ \in U_o(go)$.
		\end{cor}
		\begin{proof}
			If $go $ belongs to the interior $S_1$, then $\theta(o, go)$ is regular and by Lemma \ref{hyperbolic-creation}, $g$ is strongly regular hyperbolic. By the proof of Lemma \ref{hyperbolic-creation}, $g^{-1}o, o$ and $go$ are aligned and belong to an axis of $g$. We can then conclude by Lemma \ref{lem control attracting chamber}.
		\end{proof}

		Recall that two chambers at infinity $C$ and $C'$ are opposite if there exists a unique apartment $A$ whose boundary contain $C$ and $C'$. For any $o\in A$, we say that $C$ and $C'$ are \emph{opposite at $o$}. Clearly, two chambers are opposite if and only if there exists $o$ such that $C$ and $C'$ are opposite at $o$. Notice also that $C$ and $C'$ are opposite at $o$ if and only if their projections $\germ_o (Q(o,C)) $ and $\germ_o (Q(o,C')) $  are opposite local alcoves in the residue building $\Sigma_o X$, (this follows from \cite[Proposition 1.12]{parreau00}). We denote by $X^\infty_o (C) \subseteq \ch(\bdinf)$ the set of chambers at infinity that are opposite to $C$ at $o$.
		
		The following fact is standard, but we include a proof for completeness.
		
		\begin{lem}\label{lem opp open}
			Let $C\in \ch (\bdinf)$ and $o \in X$. The set $X^\infty_o (C) $ is open and closed.
		\end{lem} 
		
		\begin{proof}
			If $C' \in X^\infty_o (C)$, let $A$ be the unique apartment joining $C$ and $C'$, and let $x \in Q(o, C' )$ be any vertex distinct from $o$ and in the interior of $Q(o,C')$. Let $C'' \in U_o(x)$, and consider the retraction $\rho=\rho_{A, C}$ onto $A$ centered on the Weyl chamber $C$. Since $Q(x,C)\subset A$ we see that $\rho$ is the identity in restriction to $Q(x,C)$. Furthermore the image of $Q(o,C'')$ by $\rho$ is a sector of $A$ containing $o$ and $x$, and therefore is $Q(o,C')$.  It follows that the retraction $\rho$, restricted to $Q(o, C'') \cup Q(o, C)$ is an isomorphism onto $Q(o, C') \cup Q(o, C)$. Therefore $\rho (C'')=C'$ ; since retractions are distance non-increasing it follows that $C''$ and $C$ are opposite. Let $A'$ be the apartment containing $C$ and $C''$; since the retraction $\rho$ induces an isomorphism between $A'$ and $A$, and fixes $o$, we have that $o\in A\cap A'$. This proves that $C''\in X^\infty_o(C)$. Therefore $X^\infty_o(C)$ contains an open neighborhood of $C'$. We deduce that $X^\infty_o(C)$ is open.
			
			If now $C' \in \chinf $ is not opposite to $C$ at $o$, let $y$ be any vertex in the interior of $Q(o, C')$ and $U_o(y)$ is again an open set of chambers containing $C'$ that are not opposite to $C$ at $o$. Indeed, their projection on the residue building $\Sigma_o X$ is the alcove $\Sigma_o(C')$. This shows that the complement of $X^\infty_o (C)$ in $\chinf$ is open, and therefore that $X^\infty_o(C)$ is closed. 
		\end{proof}
		
		In the course of the proof, we have shown the following standard but useful fact (see also \cite[Theorem 6.2.(i)]{remy_trojan21} and its proof). 
		
		\begin{rem}\label{rem open opposite}
			Let $C,D$ be opposite chambers at infinity and let $A$ be the unique apartment joining them. Let $o \in A$ be any vertex. Let $y \in A$ be in the interior of $Q(o, C)$. Then for all chambers $C' \in \bdinf_y(D)$, $C' $ is opposite to $D$ at $o$, i.e. $\bdinf_y(D)\subseteq \bdinf_o(D)$.
		\end{rem}
		Fix $o \in X$. We denote by $\Opp_o(C) \subseteq X$ the set
		$$ \Opp_o(C):= \bigcup_{C' \in X^\infty_o (C)}  \mathrm{int}(Q(o, C')).$$
		In other words, the set $\Opp_o(C)$ is the set of all $x \in X$ that project onto an alcove that is opposite to $\Sigma_o (C) $ in the residue building $\Sigma_o X$. Because $\Opp_o(C)$  is a union of open sets, it is open for the cone topology.
		\newline 
		
		We define $\widetilde{\Opp}_o(C)$ as $\Opp_o(C) \cup X^\infty_o (C)$. It is open for the cone topology on $X \cup \bdinf$.  Lemma \ref{lem:inclusionOpp} below generalizes Remark \ref{rem open opposite}. 
		
		\begin{lem}\label{lem:inclusionOpp}
			If $y$ is a vertex in $Q(o,C)$, then $\Opp_o(C) \subseteq \Opp_y(C) $ and $ X^\infty_o (C)\subseteq X^\infty_y (C)$.
		\end{lem}
		\begin{proof}
			We write the proof for $X^\infty_o$, the case of $\Opp$ is similar, and refer the reader to Figure \ref{figure hyp creation opp} in order to follow the argument.  If $C'\in X^\infty_o(C)$, let $A$ be the unique apartment containing $C$ and $C'$. Then by definition $o\in A$ so $Q(o,C)\subset A$ and hence $y\in A$. Now $C$ and $C'$ are opposite so $C'\in X^\infty_y(C)$.

			\begin{figure}[h]
				\centering
				\begin{center}
					\begin{tikzpicture}[scale=2]
						\draw (-1,0) -- (1,0)  ;
						\draw (-1,-0.2) -- (1,-0.2)  ;
						\draw (-0.5, -0.86) -- (0.5, 0.86)  ;
						\draw (-0.7, 0.86) -- (0.3, -0.86)  ;
						\draw (-0.7, -0.86) -- (0.3, 0.86)  ;
						\draw[black] (-0.5, 0.86) -- (0.5, -0.86) ;
						\filldraw (0,0) circle(0.7pt) ;
						\filldraw (-0.32,-0.2) circle(0.7pt) ;
						\draw (0,0) node[above right]{$o$} ;
						\draw (-0.2,-0.2) node[below left]{$y$} ;
						\draw (-1.3, -0.7) node[right]{$C$} ;	
						\draw (-1.1, -0.1) to[bend right =20] (-0.75, -1) ;
						\draw [thick, draw=black, fill=gray!35]
						(1, 0) -- (0,0) -- (0.5, 0.86) ;	
						\draw (0.5, 0.4) node[right]{$\widetilde{\Opp}_o(C)$} ;
					\end{tikzpicture}
				\end{center}
				\caption{The set $\widetilde{\Opp}_o(C)$ is contained in  $\widetilde{\Opp}_y(C)$}\label{figure hyp creation opp}
			\end{figure}
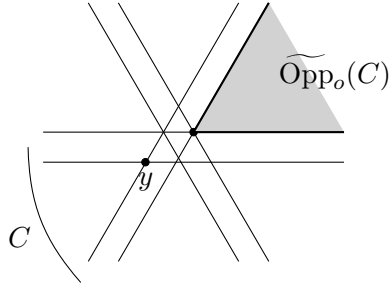
		\end{proof}	
		
		\begin{lem}\label{lem Urysohn}
			Let $o \in X$, $C \in \chinf$ and $ y $ in the interior of $Q(o, C)$. Then there exists a bounded continuous function $\varphi_{o,y,C} : X \cup \chinf \to [0,1]$ such that 
			\begin{equation}
				\varphi_{o,y,C}(x)=\begin{cases}
					1,& \text{ if } x \in \widetilde{\Opp}_o(C)\nonumber \\
					0,& \text{ if } x \notin \widetilde{\Opp}_y(C). \nonumber
				\end{cases}
			\end{equation}
		\end{lem}
		\begin{proof}

			Let the set $F$ be the closure of $\widetilde{\Opp}_o(C)$ for the cone topology, and define the set $F'$ to be $ (X\cup \bdinf) -\widetilde{\Opp}_y(C)$. Both sets are closed. By Lemma \ref{lem:inclusionOpp}, they are disjoint. The space $X \cup \bdinf $ is a normal topological space (as it is metrizable), hence Urysohn's Lemma (see for instance \cite[Theorem 4.2]{lang83}) applies and there exists a continuous function $\varphi_{o,y,C} : X \cup \chinf \to [0,1]$ satisfying 
			$$\varphi_{o,y,C}(F)= 1 \text{ and } \varphi_{o,y,C}(F')= 0,$$
			which proves the Lemma. 
		\end{proof}

		For the rest of this section, we let $\mu$ be a symmetric and admissible probability measure on $G$, and we denote by $\nu$ the unique stationary measure on $\ch(\bdinf)$ defined by Theorem \ref{thm:LB24}.  As before, we denote by $(Z_n (\omega))$ the random walk generated by $\mu$ on $G$.

		The following result is an application of \cite[Theorem 1.1]{furstenberg_kifer83}.

		\begin{lem}\label{lem:nu(omega)}
			Let $\nu$ be the unique $\mu$-stationary measure on $X \cup \chinf$. Let $C$ be a chamber at infinity and $o\in X$ be a vertex. Let  $ y $ be a point in the interior of $Q(o, C)$. Then almost surely, 	
			\begin{eqnarray}
				\lim_{n\to +\infty} \frac{1}{n} | \{ k\leq n\mid Z_k^{-1}y \in \widetilde{\Opp}_y(C) \} | \geq \nu(\widetilde{\Opp}_o(C)) \nonumber
			\end{eqnarray}
		\end{lem}
		
		\begin{proof}
			Let $\varphi_{o,y,C}$ be the function given by Lemma \ref{lem Urysohn}. As $\varphi_{o,y,C}$ is continuous and bounded,  \cite[Theorem 1.4]{furstenberg_kifer83} (with the uniqueness of the stationary measure) gives that almost surely,
			\begin{eqnarray}
				\lim_{n\to +\infty} \frac{1}{n} \sum_{k=1}^n \varphi_{o,y,C}(Z_k^{-1}y) = \nui(\varphi_{o,y,C})\nonumber.
			\end{eqnarray}
			where $\nui$ is the unique $\mui$-stationary measure on $\ch(\bdinf)$ given by Theorem \ref{thm:LB24}, and  $\mui = \iota_\ast \mu$, with $\iota(g) = g^{-1}$. 
			Since $\mu$ is assumed to be symmetric, we have $\nui=\nu$, so that the above average in fact converges to $\nu(\varphi_{o,y,C})$. 
			
			As $0\leq \varphi_{o,y,C} \leq 1$ and $\varphi_{o,y,C}(x)= 1$ if $ x\in \widetilde{\Opp}_o(C)$ we get that 
			$$\nu(\varphi_{o,y,C}) \geq \nu(\widetilde{\Opp}_o(C)).$$
			Finally, we obtain that 
			\begin{eqnarray}
				\lim_{n\to +\infty} \frac{1}{n} | \{ k\leq n\mid Z_k^{-1}y \in \widetilde{\Opp}_y(C) \} | &\geq& \lim_{n\to +\infty} \frac{1}{n} \sum_{k=1}^n \varphi_{o,y,C}(Z_k^{-1}y) \nonumber \\ &\geq& \nu(\widetilde{\Opp}_o(C)). \nonumber
			\end{eqnarray}
		\end{proof}

		\subsection{Proof of Theorem \ref{thm prop loxod}}	
		
		In this section our goal is to prove Theorem \ref{thm prop loxod}. First we will need the following lemma.
		
		\begin{lem}\label{prop opp comb conv}
			Let $o \in X$ and $C \in \ch(\bdinf)$. Let $\{x_n\}$ be a sequence of vertices in $Q(o, C)$ that converges to $C$ in the combinatorial sense. Then for all $C'  \in \ch(\bdinf)$ opposite to $C$, there exists $n_o \in \mathbb{N}$ such that for all $n \geq n_o$, $C'$ and $C$ are opposite at $x_n$. 
		\end{lem}
		
		\begin{proof}
			Let $C'$ be a chamber opposite to $C$, let $A$ be the unique apartment joining them and let $x$ be a vertex in $A$. By axiom (A4), there exists a subsector $Q(o', C)$ of both $Q(o, C)$ and $Q(x, C)$. Since $(x_n)$ converges to $C$ in the combinatorial sense and belongs to $Q(o, C)$, there exists $n_0$ such that for all $n \geq n_0$, $x_n$ belongs to the interior of $Q(o', C)$. We can then conclude by the same argument as in the proof of Lemma \ref{lem opp open}, and say that $C$ and $C' $ are opposite at $x_n $ for every such $n $.  
		\end{proof}

		\begin{rem}
			The assumption that $\{x_n\}$ remains in a sector (or in a given apartment, the proof would be the same) is really necessary. To give a counterexample, take first a sequence of vertices $(x'_n)$  in $Q(o, C)$ that converges to $C$ in the combinatorial sense. Then define $x_n$ to be a vertex at distance 1 of $x'_n$ but not in $Q(o,C)$. Then $(x_n)$ still converges to $C$ in the combinatorial sense, but if $A$ is an apartment containing $Q(o,C)$ and $C'$ is the chamber opposite to $C$ in $A$, then for every $n$ large enough the chambers $C$ and $C'$ are \emph{not} opposite at $x_n$.
		\end{rem}

		\begin{proof}[Proof of Theorem \ref{thm prop loxod}]
			First we reduce to the complete and separable case. By Proposition \ref{prop:separable} we can find
			$Y\subset X$ be a separable subbuilding which is $G$-invariant. Then the action of $G$ on $Y$ is still non-elementary, and by Proposition \ref{cor:nonelementaire} it is still non-elementary on the completion $Z=\overline Y$ of $Y$. On the other hand,   by  \cite{kleiner_leeb97} (see also \cite[Lemma 4.4]{struyve11}) we can embed $X$ into a complete Euclidean building, which we denote $\hat X$. So $Z$ is a $G$-invariant complete convex separable subspace of $\hat X$, whose boundary is a spherical building, and such that the action of $G$ is non-elementary. Hence  we can apply results of the previous section to $Z$. We write $Z^\infty$ for the boundary $\partial_\infty Z$, viewed as a spherical building. 
			For $C\in \ch(Z^\infty)$ we write $Z^\infty (C)$ for the set of chambers in the spherical building $Z^\infty$ which are opposite $C$.
			
			Let $x\in X$, and let $\Omega' \subseteq \Omega $ be the measurable set defined as all the $\omega \in \Omega$ such that $(Z_n(\omega)x)_n $ converges to a chamber $C_\omega \in \ch(Z^\infty)$. We know by Theorem \ref{thm:LB24} that $\Omega'$ can be taken to have a full measure. Let $\nu$ be the unique stationary measure on $Z \cup \ch(Z^\infty)$ given by Theorem \ref{thm:LB24}. By Theorem \ref{thm:asopposite}, two chambers at infinity are almost surely opposite. In particular, replacing again $\Omega'$ by a full measure subset, we can assume that for every $\omega \in \Omega' $, $\nu(Z^\infty (C_\omega))=1$.
			
			Fix $\omega \in \Omega'$. Take a sequence of vertices $(x_n)$ such that $x_{n+1}$ is in the interior of $Q(x_n,C_\omega)$ for every $n$, and $(x_n)$ converges to $C$ in the combinatorial sense. Applying Lemma \ref{prop opp comb conv} we see that $$Z^\infty(C_\omega) = \bigcup_{n_0\geq 0}\bigcap_{n\geq n_0} Z^\infty_{x_n}(C_\omega)$$
			
			Since $x_{n+1}$ is in the interior of $Q(x_n,C)$, we have by Lemma \ref{lem:inclusionOpp} that $Z^\infty_{x_n}(C_\omega)\subset Z^\infty_{x_{n+1}}(C_\omega)$. Therefore, $Z^\infty(C_\omega) = \bigcup_{n_0\geq 0} Z^\infty_{x_{n_0}}(C_\omega)$, and we obtain that
			$$\lim\limits_{n\to +\infty} \nu( Z^\infty_{x_n}(C_\omega)) = 1.$$
			
			Fix $\varepsilon>0$. 
			It follows that there exists an $o =o (\varepsilon,\omega) \in Z $ such that $\nu(Z^\infty_{o}(C_\omega)) \geq 1-\varepsilon$.  Fix $y$ in the interior of $Q(o, C_\omega)$. 
			Finally, fix $y' \in Z$ in the interior of $Q(y, C_\omega)$.
			By Lemma \ref{lem:nu(omega)} and Lemma \ref{lem:inclusionOpp} we have 
			\begin{eqnarray}
				\lim_{n\to +\infty} \frac{1}{n} | \{ k\leq n\mid Z_k^{-1}y' \in \widetilde{\Opp}_{y'}(C_\omega) \} | \geq \nu(\widetilde{\Opp}_{y}(C_\omega))  \geq \nu(\widetilde{\Opp}_o(C_\omega))\geq 1-\varepsilon \nonumber
			\end{eqnarray}
			By Lemma \ref{lem:inclusionOpp}, $\widetilde{\Opp}_y(C_\omega) \subseteq \widetilde{\Opp}_{y'}(C_\omega)$. By Proposition \ref{prop cv comb}, $(Z_n (\omega )y')$ converges to $C_\omega$ in the combinatorial sense: the projection of $Z_n (\omega)y'$ on $\Sigma_{y'}\hat X$ is the constant local alcove $c(\omega):= \Sigma_{y'} (C_\omega)$ for all $n$ sufficiently large, say $n_0$. For all $n\geq n_0$, Lemma \ref{hyperbolic-creation} implies that $Z_n(\omega)$ is strongly regular.
			
			It follows that
			\begin{eqnarray}
				\lim_{n\to +\infty} \frac{1}{n} | \{ k\leq n\mid Z_k \text{ is strongly regular hyperbolic} \} | \geq 1-\varepsilon \nonumber.
			\end{eqnarray}
			As $\varepsilon$ was arbitrary and $\Omega'$ is of full measure, we obtain the result. 
		\end{proof}
		
		In fact, the proof also gives a control on where the attractive/repelling fixed points of the strongly regular element lie. If $g\in \Aut(X)$ is strongly regular, then its oriented axis $\ell$ is contained in a unique apartment $A$, and its positive endpoint $\ell(+\infty)$ is contained in the interior of a unique chamber at infinity of $A$, which we call its \emph{attracting chamber} in $\ch(Z^\infty)$. Similarly $\ell(-\infty)$ is contained in a chamber at infinity called the \emph{repelling chamber} of $g$.

		Using the techniques from the proof of Theorem \ref{thm prop loxod}, we can derive the existence of a pair of independent strongly regular hyperbolic isometries in the following sense. 
		
		\begin{prop}\label{prop:controlattrep} 
			There exists a pair of strongly regular hyperbolic elements $g_1, g_2$ with attracting and repelling chambers $C_1^+, C_1^-$  and $C_2^+, C_2^-$ respectively such that $C_1^+, C_1^-,C_2^+, C_2^- \in Z^\infty$ are pairwise opposite.
		\end{prop}
		\begin{proof}
			By Theorem \ref{thm prop loxod}, there exists a strongly regular hyperbolic element $g_1$ with attracting and repelling chambers at infinity $C_1^+$ and $ C_1^-$ respectively. For every chamber $C' \in \ch(Z^\infty)$, $Z^\infty(C')$ is open in $\ch(Z^\infty)$ and by Theorem \ref{thm:asopposite}, for $\nu$-almost every $C' \in \ch(Z^\infty)$, $\nu(Z^\infty(C'))=1$. As a consequence, we can assume that $\nu(Z^\infty(C_1^+)\cap Z^\infty(C_1^-)) = 1$. By Theorem \ref{thm:LB24}, for every $U \subseteq \ch(Z^\infty)$, 
			$$ \nu(U) = \mathbb{P}(\omega \mid (Z_n(\omega)o) \text{ converges to } C', \, C' \in U),$$
			where the convergence is meant in the combinatorial sense.  
			Therefore, the measurable set
			$$\{\omega \in \Omega \mid (Z_n(\omega)o) \text{ converges to } C', \, C' \in  Z^\infty(C_1^+)\cap Z^\infty(C_1^-)\} $$
			is of full $\mathbb{P}$-measure. Call this set $\Omega'$, and let $\omega \in \Omega'$. By assumption, there exists $C_\omega \in Z^\infty(C_1^+)\cap Z^\infty(C_1^-) $ such that $(Z_n(\omega)o)$ converges to $C_\omega$. 
			Applying the same argument as before, there exists $D$ in the support of $\nu$ such that $D$ is opposite to $C_1^+$, $ C_1^-$ and $C_\omega$. 
			
			Let $A$ be the unique apartment joining $C_\omega$ and $D$, and let $o \in A$ be any vertex. Let $y \in A$ be in the interior of $Q(o, D)$. Then by Lemma \ref{lem:inclusionOpp}, the set $V :=Z^\infty_y(C_\omega)$ contains  $D$ and for all $D'\in V$, $C_\omega$ and $D'$ are opposite. Since $Z^\infty(C^+_1)\cap Z^\infty(C^-_1)$ is open, up to taking a point $y$ deeper into the Weyl chamber $Q(o, D)$, we can moreover assume that $V \subseteq Z^\infty(C^+_1)\cap Z^\infty(C^-_1)$. 
			Now since $V$ contains $D$ and $D$ belongs to the support of $\nu$, $\nu(V):=m  >0$. In particular, $$\nu(\widetilde{\Opp}_y(C_\omega)) \geq m >0.$$
			We can now apply the same argument as in the proof of Theorem \ref{thm prop loxod} to build strongly regular hyperbolic elements. By Lemma \ref{lem:nu(omega)}, 
			\begin{eqnarray}
				\lim_{n\to +\infty} \frac{1}{n} | \{ k\leq n\mid Z_k^{-1}y \in \widetilde{\Opp}_o(C_\omega) \} | \geq \nu(\widetilde{\Opp}_y(C_\omega)) \geq m \nonumber
			\end{eqnarray}
			Since $Z_n(\omega)o$ converges to $C_\omega$, by the same argument as in the proof of Theorem \ref{thm prop loxod}, we conclude that there exists an infinite number of $n $ such that $Z_n (\omega)$ is a strongly regular hyperbolic isometry. See Figure \ref{figure indep hyp} for a (1-dimensional) illustration.
			
			As moreover $C_\omega \in Z^\infty(C_1^+)\cap Z^\infty(C_1^-)\cap Z^\infty(D)$, there exists $n_0$ such that for all $n \geq n_0$, $U_o(Z_n(\omega)o) \subseteq Z^\infty(C_1^+)\cap Z^\infty(C_1^-)\cap Z^\infty(D)$. Fix $n \geq n_0$ such that $Z_n(\omega)$ is a strongly regular hyperbolic isometry and such that $Z_n^{-1}(\omega)o \in \widetilde{\Opp}_y(C_\omega)$. By Corollary \ref{cor control attr ch}, the attracting chamber at infinity $C_2^+$ of $Z_n (\omega)o$ is contained in $U_o(Z_n(\omega) o)$ and the repelling chamber at infinity $C_2^-$ belongs to $U_o(Z_n^{-1}o )$. As $y $ was arbitrarily deep inside the Weyl chamber $Q(o, D) $, we can then obtain $C_2^- \in Z^\infty(C_1^+)\cap Z^\infty(C_1^-)$ and the proposition is proven. 
		\end{proof}

		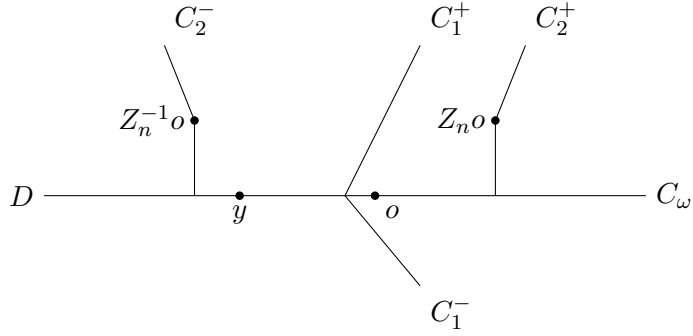
\begin{figure}[h]
			\centering
			\begin{center}
				\begin{tikzpicture}[scale=2]
					\draw (-2,0) -- (2,0)  ;
					\draw (0.5,1) -- (0,0)  ;
					\draw (0.5,-0.6) -- (0,0)  ;
					\draw (-1,0.5) -- (-1.2,1)  ;
					\draw (-1,0) -- (-1,0.5)  ;
					\filldraw (-1,0.5) circle(0.7pt) ;
					\draw (1,0.5) -- (1.2,1)  ;
					\draw (1,0) -- (1,0.5)  ;
					\filldraw (1,0.5) circle(0.7pt) ;
					\filldraw (-0.7,0) circle(0.7pt) ;
					\filldraw (0.2,0) circle(0.7pt) ;
					\draw (-0.7,0) node[below]{$y$} ;
					\draw (0.2,0) node[below right]{$o$} ;
					\draw (1,0.5) node[left]{$Z_n o$} ;
					\draw (-1,0.5) node[left]{$Z^{-1}_n o$} ;
					\draw (-2,0) node[left]{$D$} ;
					\draw (2,0) node[right]{$C_\omega$} ;
					\draw (1.2,1) node[above right]{$C_2^+$} ;
					\draw (0.5,1)  node[above right]{$C_1^+$} ;
					\draw (0.5,-0.6)  node[below right]{$C_1^-$} ;
					\draw (-1.2,1) node[above right]{$C_2^-$} ;
				\end{tikzpicture}
			\end{center}
			\caption{Construction of independent strongly regular hyperbolic isometries}\label{figure indep hyp}
		\end{figure}
		
		\subsection{Local-to-global fixed point theorem}

		First, we recall a fixed point property for trees due to Serre for simplicial trees \cite[Corollary I.6.5.3]{serre80} and to Morgan and Shalen for $\mathbb{R}$-trees \cite[Proposition II.2.15]{morgan_shalen84}.
		
		\begin{prop}\label{prop pt fixe serre}
			Let $G$ be a finitely generated group acting on a tree. If each element of $G$ has a fixed point, then $G$ has a global fixed point. 
		\end{prop}

		\begin{proof}[Proof of Theorem \ref{local-global}]
			
			By Theorem \ref{thm prop loxod}, we see that the action of $G$ on $X$ must be elementary. If $G$ has a bounded orbit then by Struyve's extension of the Bruhat-Tits fixed point theorem \cite[Main Result 1]{struyve11} we conclude that $G$ has a fixed point in $X$.
			
			If not, then $G$ has a finite orbit on $\partial_\infty X$. In this case, $G$ has a finite index subgroup $G_0$ which fixes a boundary vertex  in $\partial_\infty X$, and therefore a facet $v$ of the building at infinity. Hence $G_0$ acts on the corresponding panel tree $T_v$. The projection $\pi_v : X\to T_v$  is $G_0$-equivariant, so the action on $T_v$ is still elliptic.  By Proposition \ref{prop pt fixe serre}, $G_0$ fixes a point in $T_v$, which we denote by $\xi_0$. By \cite[Lemma 3.3]{schillewaert_struyve_thomas22} it follows that $g$ fixes some point $x^g$ in $\pi_v^{-1}(\xi_0)$. Hence it fixes a half-line $[x^g ,v)$, whose projection on $T_v$ is $\xi_0$. 
			
			Now let $\{g_1,\dots,g_n\}$ be a generating set for $G_0$. Then each $g_i$ fixes a half-line $[x^{g_i},v)$ which projects to $\xi_0$. Since all these lines have the same projection, it follows that any two of them intersect. Therefore there is a half-line contained in all the $[x^{g_i},v)$, which is fixed pointwise by all of the finitely many generators of $G_0$, and therefore by all of $G_0$. Therefore $G$ has a finite orbit on $X$, and again  by \cite[Main Result 1]{struyve11} the theorem is proved.
		\end{proof}

		\subsection{Tits alternative}
		
		The proof of Theorem \ref{thm Tits alt} is based on the following well-known "ping-pong" Lemma, which we state for cyclic subgroups.
		
		\begin{lem}\label{lem ping}
			Let $G= \langle a,b\rangle$ be a group generated by two elements acting on a set $Y$, and let $Y_1^+,Y_1^-,Y_2^+,Y_2^-$ be four disjoint non-empty subsets in $Y$, and let $Y'$ be the union of these four sets. If for every $k \geq 0$, we have
			\begin{itemize}
				\item  $a^k(Y'\setminus Y_1^-) \subseteq Y_1^+$,
				\item  $a^{-k}(Y'\setminus Y_1^+) \subseteq Y_1^-$,
				\item  $b^k(Y'\setminus Y_2^-) \subseteq Y_2^+$,
				\item $b^{-k}(Y'\setminus Y_2^+) \subseteq Y_2^-$,
			\end{itemize}		 then $G$ is isomorphic to a free group of rank 2. 
		\end{lem}
		
		The proof is similar to the original proof of Tits \cite{tits72} for linear groups. It relies on the dynamics of strongly regular hyperbolic elements. 
		
		The following lemma gives quantitative control about the speed of convergence. 
		
		\begin{lem}\label{lem quanti dyn srh}
			Let $g$ be a strongly regular hyperbolic element with translation apartment $A$ and let $C^-, C^+ \in \chinf$ be the repelling (resp. attracting) chamber at infinity for $g$. Let $o,y \in A$ be vertices such that $y \in Q(o, C^+)$. Then there exists $N$ such that for all $C \in X^\infty_o(C^-)$, $g^n (C) \in U_o(y) $ for every $n \geq N$. 
		\end{lem}
		
		\begin{proof}
			Denote by $\xi^+ \in \bd X$ the attracting fixed point of $g$. Recall that $g$ acts as a translation on $A$, and that $g^n o \to  \xi^+ $ in the conical sense. As a consequence, there exists $N$ such that for all $n \geq N$, $g^n o \in Q(y, C^+)$. Fix such $n \geq N$. For all $C \in X^\infty_o(C^-)$, the chamber $g^n (C) $ is represented by the sector $Q(g^n o, g^n C)$. As $g$ acts as a translation, and as by definition $Q(o, C^-)$ and $Q(o, C)$ are opposite sectors at $o$,  $Q(g^n o, C^-)$ and $Q(g^n o, g^n C)$ are opposite sectors at $g^n o$. As a consequence, $g^n o \in Q(y, C) $ and $g^n C \subseteq U_o(y)$. 
		\end{proof}

		\begin{proof}[Proof of Theorem \ref{thm Tits alt}]
			By assumption, the action is non-elementary. Hence, by Proposition \ref{prop:controlattrep}, there exists a pair of strongly regular hyperbolic elements $g_1, g_2 \in G$ with translation apartments $A_1,A_2$ respectively. Denote by $\xi_1^+,\xi_1^-$ (resp. $\xi_2^+,\xi_2^-$) the attracting fixed points in $\bd X$ of $g_1$ (resp. $g_2$), and by $C_{1,2}^\pm$ the chamber at infinity for which $\xi_{1,2}^\pm$ is an interior point respectively.  Moreover, 
			we can assume that the chambers $C_1^+, C_1^-,C_2^+,C_2^-$ are pairwise opposite. Let $o,o',o''\in X$. By  Lemma \ref{lem opp open}, we can choose $U_1^+$, $U_1^-,U_2^+,U_2^+$ open neighborhoods of $C_1^+, C_1^-,C_2^+,C_2^-$ respectively, such that  $U_1^+\subset X_o^\infty(C_1^-)\cap X_{o'}^\infty(C_2^+)\cap X_{o''}^\infty(C_2^-)$, and similarly for the other indices and exponents. Furthermore,  as the sets $\{U_o(y)\}_{y\in X}$ form a basis for the topology of $\ch(Z^\infty)$ for any $o$, we can choose $y_1\in Q(o,C_1^+)$, $y_2^+\in Q(o',C_1^+)$ and $y_2^-\in Q(o'',C_1^+)$ such that $U_1^+$ contains $U_o(y_1)\cup U_{o'}(y_2^+)\cup U_{o''}(y_2^-)$.
			
			Then by Lemma \ref{lem quanti dyn srh} we see that there exists $N_1^+$ such that for every $n\geq N_1^+$ we get $g_1^n(U_2^+)\cup g_1^n(U_1^+)\cup g_1^n(U_2^-)\subset U_1^+$. Arguing similarly for the other subsets, and replacing $g_1$ and $g_2$ by sufficiently large powers,  we can apply Lemma \ref{lem ping} and deduce that $g_1$ and $g_2$ generate a free group.

		\end{proof}

		\section{Elementary actions}
		
		We conclude the paper with the proof of Theorem \ref{elementary}. 
		
		\begin{proof}[Proof of Theorem \ref{elementary}]
			By Theorem \ref{thm Tits alt} we may assume the action is elementary and by Theorem \ref{semi-simple}, all elements are either elliptic or hyperbolic. If there is a global fixed point in $X$, then the group is finite since point-stabilizers are finite. So we may assume that $G$ contains a subgroup of finite index $G_0$ which fixes a vertex at infinity $v$. 
			If all elements of $G$ are elliptic then by Theorem \ref{local-global} there is a global fixed point, hence again the group is finite. Therefore we henceforth assume that there exists a hyperbolic isometry for the $G$-action on $X$. 
			
			As $G_0$ acts on the corresponding panel tree $T_v$ we get that either $G_0$ (and hence $G$) contains a non-abelian free group, or $G_0$ fixes a point $p$ of $T_v$ or an end $\xi\in \partial T_v$. 
			
			Suppose first that $G_0$ fixes a point $p\in T_v$. Let $g\in G_0$ be a hyperbolic element in the action on $X$. Let $h\in G_0$ be an elliptic element for the action on $X$. As $g$ is hyperbolic, there exists an axis $\gamma$ on which $g$ acts by translation, and since $g$ fixes $p$, we can assume that the strong asymptote class of $\gamma_{\mid [0, \infty)}$ is $p \in T_v$. Up to taking $g^{-1}$, we can assume that $v$ is the attracting fixed point of $g$ in $\bd X$. Since $h$ is elliptic,  if we denote by $o \in X$ one of its fixed points, $h$ fixes the ray $\gamma'_{\mid [0, \infty)}$ representing $[o, v)$. Now there exists $t_0, t_1$ such that $\gamma_{\mid [t_0, \infty)}$ and $\gamma'_{\mid [t_1, \infty)}$ are contained in the same apartment, and are parallel in this apartment. Since $h$ fixes $\gamma'_{\mid [t_1, \infty)}$ and also stabilizes the strong asymptote class of the ray $r:=\gamma_{\mid [t_0, \infty)}$, it must fix $\gamma_{\mid [t_0, \infty)}$ pointwise.
			Let $x$ be a point on $r$. Since $g^{-n}hg^n$ fixes $x$ for all $n\geq 0$ and since point-stabilizers are finite, there exists $k\geq 0$ such that $g^{-k}hg^k = 1$ and therefore $h=1$. Thus all non-trivial elements in $G_0$ are hyperbolic. 
			
			Let $r$ be a ray in the strong asymptote class of $p \in T_v$. Let $x \in r$ be a point. Let us denote by $S(x)$ the set 
			$$ S(x) := \{ g \in G_0 \mid x \in \mathrm{Min}(g)\}.$$ 
			Notice that for any element $g\in S(x)$, the ray $[x, v)$ belongs to some axis $\gamma_g$ of $g$, hence for $x' \in [x, v)$, $S(x) \subseteq S(x')$. Using the same argument as before, for any hyperbolic element $g \in G_0$, $g$ stabilizes a sub-ray of $r$. As a consequence, 
			$$\cup \{S(x) \mid x \in r\} = G_0.$$
			Since the action is wandering, for every $x \in r$, there exists $g_x \in S(x)$ of minimal (non-zero) translation length $\ell(g_x)$. Up to taking its inverse, we assume that $v$ is the attracting fixed point of $g_x$. 
			
			Now let $ h \in S(x)$. Again, we assume that $v$ is the attracting fixed point of $h$. Since $g_x$ is of minimal translation length, and $h$ also acts by positive translations on $[x, v)$, there exists $k \in \mathbb{N}$ such that $k\ell (g_x)= \ell(h)$. Consequently, $g_x^{-k}h(x) = x$ and since there are no non-trivial elliptic element, $h = g_x^k$. In particular, any two element $g, h \in S(x)$ commute. As $G_0$ is the increasing union of all these $S(x)$, for $x \in r$, we have proven that $G_0$ is abelian, and is $\mathbb{Z}$ or $\mathbb{Z}^2$ by the flat torus theorem \cite[Theorem II.7.1]{bridson_haefliger99}.

			We now assume that all elements $g \in G_0$ fix the end $\xi \in \partial T_v$. Due to the canonical $\mathrm{Aut}(X)_v$-equivariant bijection $\partial T_v \simeq \res(v)$ between the ends of $T_v$ and the chambers in the residue $\res(v)$ \cite[Proposition 4]{tits86}, it means that $G_0$ fixes the chamber at infinity $C:=C_\xi$ associated to $\xi$. Let us first prove that there is no non-trivial elliptic element. By contradiction, assume that there is such an element $h$. Then as $h$ fixes a point $o \in X$ and the chamber at infinity $C$, $h$ fixes the sector $Q(o, C)$. Let $g$ be a hyperbolic element in $G_0$. Let $A$ be an apartment containing $Q(o, C)$ in its boundary. By \cite[Lemma 2.2]{kramer_schillewaert17}, there exists a unique translation $t \in \iso(A)$ of $A$ such that for all $p \in gA \cap A $, $g^{-1}(p )= t(p)$. As $g^{-1}$ also fixes $C$, there also exists $s \in \iso(A)$ a unique translation of $A$ such that for all $p \in A \cap g^{-1}A$, $g(p) = s(p)$. In particular, for all $p \in  A \cap g^{-1}A$, $g^{-1}\circ g(p) = t\circ s(p) = p $, therefore $s = t^{-1}$. Then, taking $p$ sufficiently deep in $Q(o, C)$ (for instance, $p \in gA \cap A \cap t^{-1}(Q(o, C))$), we have that $ghg^{-1}(p) = p$. Since pointwise-stabilizers are uniformly bounded, there exists $q$ such that $g^q hg^{-q} = 1$ and therefore $h=1$.
			
			Now let $h \in G_0 $ be another element. We denote by $s$ the translation of $A$ such that for all $p \in A \cap hA$, $h^{-1}(p)=s(p)$ as before. For $p $ sufficiently deep in any Weyl chamber $S$ representing $C$, we obtain that $ghg^{-1}h^{-1}(p) = tst^{-1}s^{-1}(p) = p$ (in fact, this holds for any $p \in A \cap hA \cap hgA \cap hgh^{-1}A \cap hgh^{-1}g^{-1}A$, which is non-empty because all elements in $G_0$ fix $C$). As there are no non-trivial elliptic elements, we have that $[g,h]=1$. Then $G_0$ is abelian, and is $\mathbb{Z}$ or $\mathbb{Z}^2$ by the flat torus theorem \cite[Theorem II.7.1]{bridson_haefliger99}.
		\end{proof}

		\bibliographystyle{alpha}
		
		\bibliography{bibliography}

	\end{document}